\documentclass[12pt]{amsart}
\usepackage{amsmath,amsfonts,amssymb,amscd,amsthm,amsbsy,epsf}
\newtheorem*{THM}{Theorem}
\newtheorem*{COR}{Corollary}

\def\qed{\hbox{\hskip 6pt\vrule width6pt height7pt depth1pt \hskip1pt}}

\def\Z{{\mathbb Z}}
\def\H{{\bf H}}
\def\R{{\mathbb R}}

\begin{document}

\title{Tiling Spaces are Inverse Limits}
\author{Lorenzo Sadun}
\begin{abstract}
Let $M$ be an arbitrary Riemannian homogeneous space, and let $\Omega$ be
a space of tilings of $M$, with finite local complexity (relative to some
symmetry group $\Gamma$) and closed in the natural topology.  Then 
$\Omega$ is the inverse limit of a sequence of compact finite-dimensional
branched manifolds.  The branched manifolds are (finite) unions of cells, 
constructed from the tiles themselves and the group $\Gamma$. 
This result extends previous results of Anderson and
Putnam \cite{ap}, of Ormes, Radin and Sadun \cite{ors}, of Bellissard,
Benedetti and Gambaudo \cite{bbg}, and of G\"ahler \cite{gaehler}. 
In particular, the construction in this paper 
is a natural generalization of G\"ahler's.
\end{abstract}
\address{Department of Mathematics, The University of Texas at Austin,
Austin, TX 78712}
\email{sadun@math.utexas.edu}
\subjclass{Primary: 37B50, Secondary: 52C23, 54F65, 37C85}
\maketitle


\markboth{Lorenzo Sadun}{Tiling Spaces are Inverse Limits}


\section{Background}

In the last few years, it has become clear that many spaces of tilings
of $\R^d$ can be viewed as inverse limit spaces.  Anderson and Putnam
\cite{ap} began this program for substitution tilings. Given a
substitution, they showed that the corresponding space of tilings of
$\R^d$ is the inverse limit of a branched $d$-manifold $K$ under an
expansive map from $K$ to itself.  If the substitution has a
property called ``forcing the border'' \cite{kellendonk}, then the
manifold $K$ is constructed by stitching all the tile types together
along possible common boundaries.  If the substitution does not force
the border, then the construction is similar, only using collared
tiles.  (A collared tile is a tile that is labeled by the pattern of
tiles that touch it). 
For this construction to work, the tilings must involve only a
finite number of tile types (up to translation), meeting full-face to
full-face.  In particular, the construction does not apply to tilings
like the pinwheel \cite{pinwheel}, where tiles appear in an infinite
number of orientations.

Ormes, Radin and Sadun \cite{ors} extended the Anderson-Putnam 
construction to substitution tilings of $\R^d$ on which the entire
Euclidean group acts continuously.  Tiles may appear in arbitrary
orientations, but there can only be a finite number of tile types
{\em up to Euclidean motion}, and tiles must meet full-face to full-face. 
The branched manifold has dimension $d(d+1)/2$, which is the dimension
of the Euclidean group.  

In this construction, a cell in the branched manifold $K$ is not a
tile.  Rather, a cell is the product of a (possibly collared) tile
with $SO(d)$, modulo any (finite!) rotational symmetry that the tile
might have.  This gives a description of all the ways a tile
containing the origin may be placed.  The substitution (call it
$\sigma$) replaces each oriented tile with a union of oriented tiles,
giving a map from $K$ to itself.
Such a union of tiles is called a {\em supertile} of order 1.  The
substitution applied to a supertile of order 1 gives a supertile of
order 2, and so on.  

A point $(x_0, x_1, \ldots)$ in the inverse limit $\leftarrow_\sigma
K$ is a consistent description of a tiling, with $x_0$ telling how the
origin sits inside a tile, $x_1$ telling how the origin sits inside a
supertile of order 1, and $x_n$ telling how the origin sits inside a
supertile of order $n$.  If the substitution forces its border (or if
we are using collared tiles), the sequence $(x_0, x_1, \ldots)$ gives
a consistent description of a unique tiling of $\R^d$.

More recently, G\"ahler \cite{gaehler} and Bellissard, Benedetti and
Gambaudo \cite{bbg} have each applied inverse limit methods to tilings
that need not be generated by a substitution.  If $x$ is a tiling of
$\R^d$ that has finitely many tile types {\em up to translation},
meeting full-face to full-face, then the continuous hull of $x$ (i.e.,
the closure of the translational orbit of $x$) is the inverse limit of
a sequence of compact branched manifolds $K_0$, $K_1$, $K_2$,
$\ldots$, under a sequence of maps $\sigma_n: K_n \to K_{n-1}$, where
each branched manifold $K_n$ is the union of (marked) tiles from the
original tiling.  Of the two constructions, G\"ahler's is conceptually
simpler, but that of Bellissard, Benedetti and Gambaudo appears to be
calculationally stronger, leading to results such as gap-labeling
theorems \cite{bbg}.

This paper is an extension of G\"ahler's construction to tilings of arbitrary
Riemannian homogeneous spaces, with general symmetry group.  The generalization
of the Bellissard-Benedetti-Gambaudo approach to arbitrary spaces is 
being done independently by Benedetti and Gambaudo \cite{gambaudo}.

\section{Theorem and Proof}

Before stating and proving the result, we must establish some notation.
Let $M$ be a Riemannian homogeneous space (such as $\Z^d$, $\R^d$, $\H^2$, 
$\H^2 \times \R^3$, etc.), and pick a point to be the origin.
Let $G$ be the group of isometries of $M$, let $\Gamma$ be a 
closed subgroup of $G$, and let $\Gamma_0$ be the
subgroup of $\Gamma$ that fixes the origin.  Let $\Omega$ be a 
collection of tilings of $M$. We give $\Omega$ the topology that two
tilings are $\epsilon$-close if they agree on a ball of size $1/\epsilon$
around the origin, up to the action of an $\epsilon$-small 
element of $\Gamma$.  We assume that $\Omega$ is closed under the action
of $\Gamma$ (i.e., $\Omega$ is a union of $\Gamma$-orbits), and that
$\Omega$ is compact.  This implies that $\Omega$ has finite local complexity,
up to the action of $\Gamma$. 

\begin{THM} \label{mail} $\Omega$ is the inverse limit 
of a sequence of compact branched manifolds $K_1, K_2, \ldots$ and
continuous maps $\sigma_n: K_n \to K_{n-1}$. The dimension of the branched
manifold is the dimension of $\Gamma$. 
\end{THM}

The idea of the proof is quite simple.  A point in the $n$-th
approximant $K_n$ is a description of a tile containing the origin,
its nearest neighbors (sometimes called the ``first corona''), its
second nearest neighbors (the ``second corona'') and so on out to the
$n$-th nearest neighbors.  (For these purposes, tiles that meet at a point
are considered nearest neighbors.) The map $\sigma_n: K_n \to K_{n-1}$ simply
forgets the $n$-th corona.  A point in the inverse limit is then a
consistent prescription for constructing a tiling out to infinity.  In
other words, it is a tiling.

What remains is to actually construct $K_n$ out of geometric pieces and
show that $K_n$ is a branched manifold.  

First suppose that the tiles are polytopes that meet
full-face to full-face.  We consider two tiles $t_1$, $t_2$ in (possibly
different) tilings of $M$ to be
equivalent if a patch of the first tiling, containing $t_1$ and its 
first $n$ coronas, is identical, up the the action of $\Gamma$, to a 
similar patch around $t_2$.   Since $\Omega$
has finite local complexity, there are only finitely many
equivalence classes, each of which is called an $n$-collared tile. 

For each $n$-collared tile $t_i$, we consider how such a tile
can be placed around the origin.  Let $s_i\subset t_i$ be the set of  
points where the origin may sit.  By finite local complexity, there can
only be a finite number of connected components to $s_i$, and each component
is a submanifold of $t_i$ with the same dimension as $\Gamma/\Gamma_0$. 
If $t_i$ does not admit any symmetry, then for each point $p \in s_i$,
$\Gamma_0$ acts simply transitively on the ways to place $t_i$ down with
the spot $p$ landing at the origin.  The set of ways to place $t_i$ is
therefore a principal $\Gamma_0$ bundle over $s_i$, which we denote $E_i$.  
The cell $C_i \subset K_n$ associated with $t_i$ is then exactly
$E_i$.  

If there are no 
topological obstructions to trivializing this bundle, 
we make the identification
\begin{equation} C_i = E_i = s_i \times \Gamma_0. \end{equation}
If $M$ is flat, then there is a {\em canonical} trivialization of the
frame bundle, and this descends to a canonical product (1). If $\Gamma$
acts transitively on $M$, then $s_i = t_i$ is contractible, and 
the decomposition (1), while not canonical, is guaranteed to exist.  
Although there do exist tilings
where neither of these conditions are met, the author knows of no 
examples where $C_i$ fails to be trivializable. 

If $t_i$ admits a discrete symmetry (e.g., is a regular $n$-gon in a 
tiling of $\R^2$ or $\H^2$), then more than one point in $s_i \times \Gamma_0$
may describe the same placement of a tile containing the origin.  In that
case, the cell associated to $t_i$ is the quotient of the $\Gamma_0$ bundle
$E_i$ by the symmetry.  That is, 
\begin{equation} C_i 
= E_i/\Gamma_{t_i} \quad (\ = s_i \times_{\Gamma_{t_i}} \Gamma_0, \hbox{ if $E_i$ 
is trivializable}),
\end{equation} where $\Gamma_{t_i} \subset \Gamma_0$ is the group of
symmetries of $t_i$.  Since $t_i$ is a collared tile, $\Gamma_{t_i}$
must be a discrete subgroup of $\Gamma_0$. (Even if a tile had a
continuous symmetry, its first corona could not.)  By construction,
$\Gamma_{t_i}$ acts without fixed points on $E_i$, so the interior of
$C_i$ is indeed a manifold.  (For instance, if $M=\R^2$ and $\Gamma$
is the Euclidean group, then $C_i$ is a Seifert fibered space.  
There may be multiple fibers over points of symmetry, 
but the total space is smooth.)

A patch of a tiling in which the origin is on the boundary of two or more tiles
is described by points on the boundary of two or more cells, 
and these points must be identified.  The branched manifold $K_n$ is 
the disjoint union of the cells $C_i$, modulo this identification. 
Since we are using $n$-collared tiles with $n\ge 1$, 
each of the points being identified
carries complete information about the placement of all the tiles that meet
the origin, together with their first $n-1$ coronas. 

We must show that a neighborhood of such a branch point is the
union of topological disks whose tangent spaces may be identified.  Each
such disk is obtained by taking a patch of a tiling in which the above data
is actually realized, and considering its orbit under the action of a
neighborhood of the identity in $\Gamma$.  This shows that the dimension of
$K_n$ is the dimension of $\Gamma$. 

Finally, we remove the assumption that the tiles are polytopes that
meet full-face to full-face.  To a tiling by other shapes we may
associate a pattern of marked points, where a special point is chosen
from each tile and labeled by the type of that tile.  The Voronoi
cells of those points are then polytopes whose faces, properly
subdivided, meet full-face to full-face.  The original tiling and the 
tiling by Voronoi cells are mutually locally derivable \cite{BSJ}, and so are
described by the same topological space, and hence by the same inverse
limit structure. \qed

In this construction, the group $\Gamma_0$ acts naturally on each
space $K_n$, and the maps $\sigma_n$ are equivariant, from which it
follows that

\begin{COR} \label{orbifolds} The space $\Omega/\Gamma_0$ of tilings
modulo rotation is the inverse limit of a sequence of compact branched
orbifolds $K_n/\Gamma_0$. 
\end{COR}

\section{Examples}

\begin{enumerate}
\item {} If $M=\Gamma=\Z^d$, then we have a $\Z^d$ subshift. The total
space $\Omega$ is a Cantor set.  The $n$-th approximant $K_n$ is a
finite collection of points, corresponding to a decomposition of the
Cantor set into a finite number of clopen sets.  This decomposition
becomes finer as $n \to \infty$, and the Cantor set is recovered as
the inverse limit. 

\item {} If $M = \R^d$ and $\Gamma=\Z^d$, then (up to a fixed translation)
$\Omega$ is  a space of tilings
of $\R^d$ by square tiles centered at the lattice points.  This is a different
description of the previous example.   
In these examples, note that $\Omega$ does not have to be the hull
of a single tiling, and that the $\Z^d$ action need not be minimal. The 
$\Z$ subshift on two letters, in which one of the letters appears at 
most twice, is neither minimal nor the closure of a single orbit, but
is an inverse limit space.

\item {} The ($d$-fold) suspension of a $\Z^d$ subshift has
$M=\Gamma=\R^d$.  This is a space of tilings of $\R^d$ by unit cubes
oriented parallel to the coordinate axes.  

\item {} A $\Z^d$ subshift may be suspended in some directions but not in
others.  For instance, the suspension of a $\Z^2$ subshift in the $x$ 
direction is a space of tilings of $\R^2$ by square tiles, meeting
full-face to full-face, whose centers have integral $y$ coordinate.
In this case $M=\R^2$ and $\Gamma = \R \times \Z$.  

\item {} The Penrose tiling space, or any other tiling of $\R^d$ with
a finite set of prototiles up to translation, has $M=\Gamma=\R^d$.
Since $\Gamma_0$ is trivial and $\Gamma$ is the full translation
group, the cells $C_i$ can be identified with the collared tiles $t_i$
themselves. This is precisely G\"ahler's construction. As was shown in
\cite{sadun-williams}, such a space is homeomorphic to the suspension
of a $\Z^d$ subshift.

\item {} The pinwheel tiling space \cite{pinwheel} has $M=\R^2$ and
$\Gamma$ the 2-dimensional Euclidean group \cite{ors}. 

\item {} In tiling hyperbolic space, there are a number of interesting 
choices for $\Gamma$.  If $\Gamma$ is a discrete group, then we have the
analog of a subshift, associating letters to a discrete set of points 
in the space being tiled.  At the other extreme, one can take $\Gamma$
to be the entire group of isometries of $\H^n$.

\item{} One dimensional orientable hyperbolic attractors are 
either solenoids or one
dimensional tiling spaces \cite{will,ap}.  
However, the dyadic solonoid {\em can} be 
viewed as a tiling space, of $\H^2$ rather than $\R^1$, 
following a construction of Penrose \cite{penrose-hyperbolic}.  See
figure 1. In the upper-half-plane model, the basic tile looks like a
rectangle, with the sides of the rectangle geodesics, 
with the top and
bottom edges horocyclic, and with the size chosen such that
the bottom edge has twice the length of the top edge.  Here the group is
$\Gamma=\Z \ltimes \R$, acting on $\H^2$ by $(n,t) (x,y) = (t+ 2^n x,2^n y)$.

\begin{figure}
\vbox{
\centerline{\epsfysize=10.0cm \epsfbox{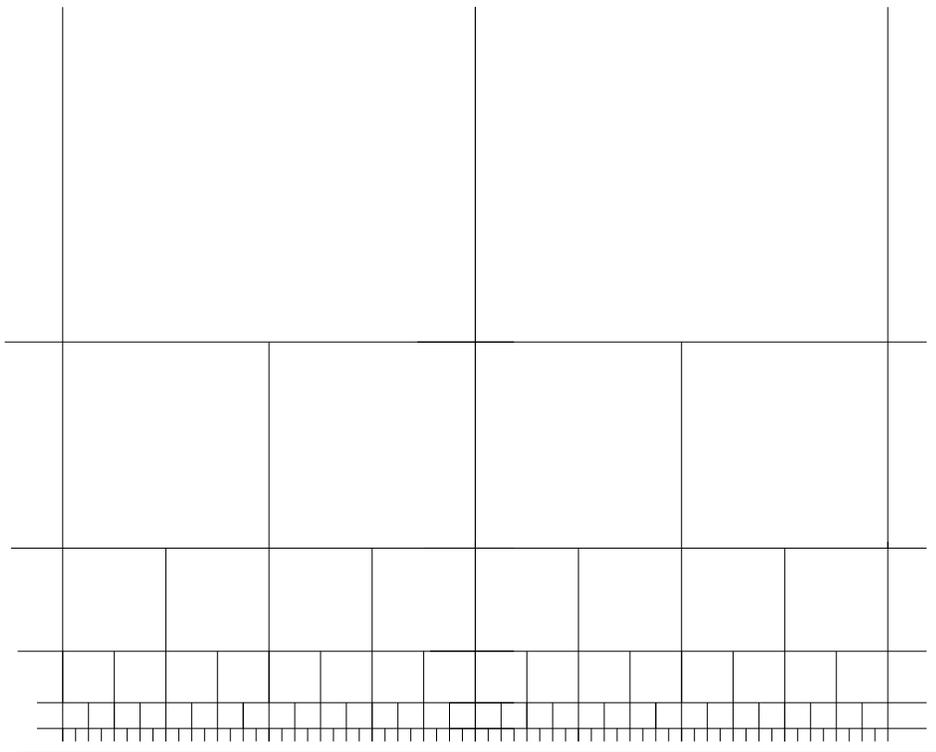} }}
\caption{Penrose's dyadic tiling of hyperbolic space}
\end{figure}

\item {} More generally, any geometric substitution in $\R^d$ gives rise
to a space of tilings of $\H^{d+1}$, with group $\Gamma =
\Z \ltimes \R^d$.  As with the dyadic solenoid, it doesn't matter whether
the substitution is invertible, since the $\Z$ action enforces the
hierarchy. Chaim Goodman-Strauss has adapted this construction to produce
a strongly aperiodic set of prototiles for $\H^2$ \cite{gs1}, and to
develop a general formalism for describing tilings of hyperbolic space
\cite{gs2}.  
\end{enumerate}

\section{Conclusions and open problems}

The inverse limit structure of $\Omega$ implies that the \v Cech
cohomology $H^*(\Omega)$ is the direct limit of $H^*(K_n)$ under the
pullback maps $\sigma_n^*$.  Every element of $H^*(\Omega)$ is the
pullback, under the natural projection $\pi_n: \Omega \to K_n$, of a
cohomology class in $K_n$, for $n$ sufficiently large.  If (and only
if) $H^*(\Omega)$ is finitely generated, then for $n$ large enough the
entire cohomology of $\Omega$ is the quotient of $H^*(K_n)$ by the
kernel of $\pi_n^*$.

To make effective use of this principle, however, requires specific
knowledge of the tiling space in question.  For substitution tilings,
it is easiest to work with the Anderson-Putnam inverse limit
construction, rather than that constructed here, although in fact the
two are shift equivalent.  For cut-and-project tilings with
sufficiently nice ``windows'', G\"ahler \cite{gaehler} has shown that
$\pi_n^*$ is actually an isomorphism in cohomology for $n$
sufficiently large, with the required size of $n$ computable from the
geometry of the window.

The inverse limit structure of tiling spaces is related to a possible
fiber bundle structure.  Locally, $\Omega$ looks like a piece of
$\Gamma$ times a Cantor set.  Can these neighborhoods be stitched
together to yield a fiber bundle (with Cantor set fiber) over a
compact manifold?  Is that manifold the quotient of the identity
component of $\Gamma$ by a co-compact subgroup? When $M=\Gamma=\R^d$,
the answer to both questions is yes \cite{sadun-williams}, but the
general case is not known.

\section{Acknowledgements}

Many thanks to the organizers and participants in the August 2002
conference on Aperiodic Order, Dynamical Systems, Operator Algebras
and Topology at the University of Victoria, British Columbia.  The
results of this paper are a direct consequence of the talks at that
conference (especially those of Franz G\"ahler and Jean-Marc Gambaudo)
and the discussions that ensued.  Additional thanks to Charles Radin
and Michael Baake for assistance with the manuscript.  This work was
partially supported by the Texas Advanced Research Program.

\end{document}